\documentclass{amsart}
\usepackage{amssymb}
\usepackage{amsmath}
\usepackage{a4}
\begin{document}
\def\id{\operatorname{id}}
\newtheorem{theorem}{Theorem}[section]
\newtheorem{lemma}[theorem]{Lemma}
\newtheorem{remark}[theorem]{Remark}
\newtheorem{definition}[theorem]{Definition}
\newtheorem{corollary}[theorem]{Corollary}
\newtheorem{example}[theorem]{Example}
\def\qedbox{\hbox{$\rlap{$\sqcap$}\sqcup$}}
\makeatletter
  \renewcommand{\theequation}{%
   \thesection.\alph{equation}}
  \@addtoreset{equation}{section}
 \makeatother
\def\BB{\mathcal{B}}
\def\MM{{\mathfrak{M}}}
\title[Pseudo-Riemannian manifolds with commuting curvature operators]
{Pseudo-Riemannian manifolds with Commuting Jacobi Operators}
\author{M. Brozos-V\'azquez and P. Gilkey}
\begin{address}{MBY:Department of Geometry and Topology, Faculty of Mathematics, University of
Santiago de Compostela, 15782 Santiago de Compostela, Spain}\end{address}
\begin{email}{mbrozos@usc.es}\end{email}
\begin{address}{PG: Mathematics Department, University of Oregon, Eugene, OR 97403, USA}\end{address}
\begin{email}{gilkey@uoregon.edu}\end{email}
\begin{abstract} We study the geometry of pseudo-Riemannian manifolds which are
Jacobi--Tsankov, i.e.
${\mathcal{J}}(x){\mathcal{J}}(y)={\mathcal{J}}(y){\mathcal{J}}(x)$
for all $x,y$. We also study manifolds which are $2$-step Jacobi
nilpotent, i.e. ${\mathcal{J}}(x){\mathcal{J}}(y)=0$ for all
$x,y$.
\end{abstract}
\maketitle

\section{Introduction}\label{sect-1}

Let $\mathcal{M}:=(M,g)$ be a pseudo-Riemannian manifold of
signature $(p,q)$ and dimension $m=p+q\ge3$; $\mathcal{M}$ is said
to be {\it Riemannian} if $p=0$ and {\it Lorentzian} if $p=1$.
Although the Riemannian and Lorentzian settings are perhaps the
most frequently studied, pseudo-Riemannian manifolds with other
signatures are important in many physical applications;
see, for example, the discussion of Kaluza-Klein gravity in
Overduin and Wesson \cite {OW97} or the brane world cosmology of
Shtanov and Sahni \cite{SS02}. Thus the higher signature setting
is important not only mathematically, but also in physical
applications.

Let
$\mathcal{R}$ be the curvature operator and $\mathcal{J}$ the Jacobi operator which are defined by the Levi-Civita
connection on $\mathcal{M}$:
\begin{eqnarray*}
&&\mathcal{R}(x,y):=\nabla_x\nabla_y-\nabla_y\nabla_x-\nabla_{[x,y]},\\
&&\mathcal{J}(x):y\rightarrow\mathcal{R}(y,x)x\,.
\end{eqnarray*}

The relationship between the spectral geometry of $\mathcal{J}$ and
the underlying geometry of the manifold has been studied extensively in recent years.
Suppose that
$\mathcal{M}$ is Riemannian. If
$\mathcal{M}$ is a
$2$-point homogeneous space, then the group of isometries acts transitively on the unit sphere bundle
$S(\mathcal{M})$ and hence the eigenvalues of $\mathcal{J}$ are constant on $S(\mathcal{M})$. Osserman \cite{Os90}
wondered if the converse is true, at least locally. He conjectured that if $\mathcal{M}$ is a Riemannian manifold
such that the eigenvalues of
$\mathcal{J}$ are constant on $S(\mathcal{M})$, then either $\mathcal{M}$ is flat or $\mathcal{M}$ is locally
isometric to a rank $1$-symmetric space. This conjecture has been established in dimensions $m\ne16$ by the work
of Chi
\cite{Chi88} and Nikolayevsky \cite{N04,N05}; the case $m=16$ is still open.

Let $S^\pm(\mathcal{M})$ be the pseudo-sphere bundles of unit
spacelike ($+$) or unit timelike ($-$) vectors. One says that a
pseudo-Riemannian manifold $\mathcal{M}$ is {\it spacelike
Osserman} (resp. {\it timelike Osserman}) if the eigenvalues of
the Jacobi operator $\mathcal{J}$ are constant on
$S^+(\mathcal{M})$ (resp. on $S^-(\mathcal{M})$). Work of
Garc\'{\i}a--R\'{\i}o et. al. \cite{GKV97} shows these are
equivalent concepts so one simply speaks of an Osserman manifold.
It is known \cite{BBG97,GKV97} that any Lorentzian Osserman
manifold has constant sectional curvature; thus the geometry is
very rigid in this setting. However if $p\ge2$ and $q\ge2$, there
are Osserman pseudo-Riemannian manifolds which are not locally
homogeneous; see, for example, \cite{BBGZ97,DGV05}.

One can weaken this condition slightly. Let $p\ge1$ and $q\ge1$.
One says that $\mathcal{M}$ is {\it pointwise Osserman} if the
spectrum of $\mathcal{J}$ is constant on $S_P^+(\mathcal{M})$, or
equivalently on $S^-_P(\mathcal{M})$, for every $P\in M$. Bla\v
zi\'c \cite{B05} has shown that if the spectrum of $\mathcal{J}$
is bounded on either $S^+_P(\mathcal{M})$ or, equivalently,
$S^-_P(\mathcal{M})$, for every $P\in M$, then necessarily
$\mathcal{M}$ is {\it pointwise Osserman}.

In this paper, instead of focusing on the spectrum, we will relate
commutativity properties of $\mathcal{J}$ to the underlying
geometry.
\begin{definition}\label{defn-1}\rm One says that a pseudo-Riemannian manifold $\mathcal{M}$ is:
\begin{enumerate}
\item {\it $2$-step Jacobi nilpotent} if $\mathcal{J}(x)\mathcal{J}(y)=0$ for all tangent vectors $x,y$.
\item {\it Jacobi--Tsankov} if
$\mathcal{J}(x)\mathcal{J}(y)=\mathcal{J}(y)\mathcal{J}(x)$ for all tangent vectors $x,y$.\item
 {\it Orthogonally Jacobi--Tsankov} if
$\mathcal{J}(x)\mathcal{J}(y)=\mathcal{J}(y)\mathcal{J}(x)$ for all $x\perp y$.
\end{enumerate}\end{definition}

Clearly (1) $\Rightarrow$ (2) $\Rightarrow$ (3). The following seminal result was established by Tsankov \cite{T05}:

\begin{theorem}\label{thm-1.2}
Let $\{\lambda_i\}$ be the eigenvalues of the shape operator of a hypersurface $M$ in $R^{m+1}$.
Then $M$ is orthogonally Jacobi--Tsankov if and only if either $\lambda_1 = ...= \lambda_m$ or $\lambda_1=
...=\lambda_{m-1}=0, \ \lambda_mΠ\ne 0$.
\end{theorem}

Theorem \ref{thm-1.2} has been extended from hypersurfaces to the more general setting in \cite{BG05}:

\begin{theorem}\label{thm-1.3}
Let $\mathcal{M}$ be an orthogonally Jacobi--Tsankov Riemannian manifold.
Then $\mathcal{M}$ has constant sectional curvature.
\end{theorem}

In passing to more general signatures, we shall impose a stronger
condition and study Jacobi--Tsankov manifolds. It is convenient to
work in the algebraic context. Let $V$ be a finite dimensional
real vector space. Let $\mathfrak{A}(V)\subset\otimes^4V^*$ be the
space of {\it algebraic curvature tensors}; these are the
$4$-tensors with the same symmetries as the Riemann
curvature tensor. Thus $A\in\mathfrak{A}(V)$ if and only if we
have the following symmetries for all $x,y,z,w\in V$:
\begin{eqnarray*}
&&A(x,y,z,w)=-A(y,x,z,w)=A(z,w,x,y),\\
&&A(x,y,z,w)+A(y,z,x,w)+A(z,x,y,w)=0\,.
\end{eqnarray*}

Let $\mathfrak{M}:=(V,\langle\cdot,\cdot\rangle,A)$ where $A\in\mathfrak{A}(V)$ and where
$\langle\cdot,\cdot\rangle$ is a non-degenerate symmetric bilinear form of signature
$(p,q)$ on $V$ which is used to raise and lower indices. The corresponding {\it algebraic curvature operator}
$\mathcal{A}\in V^*\otimes V^*\otimes\operatorname{End}(V)$ is characterized by
$$\langle\mathcal{A}(x,y)z,w\rangle=A(x,y,z,w)$$
and the Jacobi operator $\mathcal{J}=\mathcal{J}_A$ is given by $\mathcal{J}(x):y\rightarrow\mathcal{A}(y,x)x$. The
notions of Definition \ref{defn-1} then extend to the algebraic setting.
In Section \ref{sect-2}, we will show that:

\begin{theorem}\label{thm-1.4}
Let $\mathfrak{M}$ be Jacobi--Tsankov. Then:
\begin{enumerate}
\item ${\mathcal{J}}(x)^2=0$ for all $x\in V$. \item
$\mathfrak{M}$ is Osserman. \item If $V$ is Riemannian or
Lorentzian, then $A=0$.
\end{enumerate}
\end{theorem}

We can draw the following geometrical consequence from Theorem \ref{thm-1.4}:

\begin{corollary}\label{cor-1.5}
Let $\mathcal{M}$ be a Jacobi--Tsankov pseudo-Riemannian manifold of signature $(p,q)$. Then $\mathcal{M}$ is
nilpotent Osserman. If $p=0$ or if $p=1$, then $\mathcal{M}$ is flat.
\end{corollary}

One might conjecture that the condition $\mathcal{J}(x)^2=0$ for
all $x\in V$ is sufficient to imply $\mathfrak{M}$ is
Jacobi--Tsankov. This is in fact not the case as we will show in
Lemma \ref{lem-2.2}.

It is clear that any $2$-step Jacobi nilpotent algebraic curvature tensor is Jacobi--Tsankov. In Section \ref{sect-3}, we will show
that the converse holds in low dimensions:
\begin{theorem}\label{thm-1.6}
Let $\mathfrak{M}$ be Jacobi--Tsankov. If $\dim(V)\le13$, then $\mathfrak{M}$ is $2$-step Jacobi nilpotent.
\end{theorem}

The condition ${\dim}(V)\leq 13$ in Theorem \ref{thm-1.6} is
sharp. In Lemma \ref{lem-3.2}, we construct a Jacobi--Tsankov
tensor in signature $(8,6)$, which is indecomposable and for which
there exist $(x,y)$ so that $\mathcal{J}(x)\mathcal{J}(y)\neq 0$.

There are similar questions for the skew-symmetric curvature
operator.

\begin{definition}\label{defn-7}\rm
One says that ${\mathfrak{M}}$ is:
\begin{enumerate}
\item {\it $2$-step skew-curvature nilpotent} if
${\mathcal{A}}(x_1,x_2){\mathcal{A}}(x_3,x_4)=0$ for all tangent
vectors $x_1,x_2,x_3,x_4$. \item {\it Skew--Tsankov} if
${\mathcal{A}}(x_1,x_2){\mathcal{A}}(x_3,x_4)={\mathcal{A}}(x_3,x_4){\mathcal{A}}(x_1,x_2)$
for all tangent vectors $x_1,x_2,x_3,x_4$.
\end{enumerate}
\end{definition}

Motivated by Theorem \ref{thm-1.6},  in Section \ref{sect-4}, we
will study $2$-step Jacobi nilpotent algebraic curvature tensors
in relation to $2$-step skew-curvature nilpotent ones. If
$A_W\in\mathfrak{A}(W)$, we say that $(W,A_W)$ is {\it
indecomposable} if there is no decomposition
$(W,A_W)=(W_1,A_1)\oplus (W_2,A_2)$ where $\dim(W_i)\ge1$.
Similarly, we say that $\mathfrak{M}$ is indecomposable if there
is no decomposition
$\mathfrak{M}=\mathfrak{M}_1\oplus\mathfrak{M}_2$ so that
$\dim(V_i)\ge1$.

\begin{definition}\label{defn-1.8}\rm
Let $A_W\in\mathfrak{A}(W)$. Assume that $(W,A_W)$ is
indecomposable. Let $\{\bar e_1,...,\bar e_k\}$ be a basis for an
auxiliary vector space $\bar W$. Let
\begin{equation}\label{eqn-1.a}
\begin{array}{l}
\mathfrak{M}:=(W\oplus\bar W,
\langle\cdot,\cdot\rangle_{W\oplus\bar
W},A_W\oplus 0)\quad\text{where}\\
\langle e_i,e_j\rangle=\langle\bar e_i,\bar e_j\rangle=0,\quad\langle e_i,\bar
e_j\rangle=\delta_{ij}\,.\vphantom{\vrule height 12pt}
\end{array}\end{equation}
\end{definition}

We will establish the following classification theorem:

\begin{theorem}\label{thm-1.9} The following statements are
equivalent:
\begin{enumerate}
\item $\mathfrak{M}$ is $2$-step Jacobi nilpotent and
indecomposable, \item $\mathfrak{M}$ is $2$-step skew-curvature
nilpotent and indecomposable, \item $\mathfrak{M}$ is isomorphic
to one of the tensors described in Definition \ref{defn-1.8}.
\end{enumerate}
\end{theorem}

One has the following geometrical examples which arose in the study of Osserman manifolds. We refer to \cite{GIZ02,GN04} for further
details.
\begin{theorem}\label{thm-1.10} Let
$(x_1,...,x_p,y_1,...,y_p)$ be coordinates on $\mathbb{R}^{2p}$ for $p\ge2$. Let $\psi_{ij}(x)=\psi_{ji}(x)$ be a
symmetric
$2$-tensor. Let
\begin{eqnarray*}
&&g_\psi(\partial_{x_i},\partial_{x_j})=\psi_{ij}(x),\quad
g_\psi(\partial_{x_i},\partial_{y_j})=\delta_{ij},\quad
g_\psi(\partial_{y_i},\partial_{y_j})=0\,.
\end{eqnarray*}
Then $\mathcal{M}:=(\mathbb{R}^{2p},g_\psi)$ is a complete
pseudo-Riemannian manifold of neutral signature $(p,p)$ which is
$2$-step Jacobi nilpotent {  and $2$-step skew-curvature
nilpotent.}
\end{theorem}

\section{The proof of Theorem \ref{thm-1.4}}\label{sect-2}

The Jacobi operator is quadratic in $x$. We polarize to define an operator valued bilinear form by
setting:
$$
\mathcal{J}(x,y):z\rightarrow
 \textstyle\frac12\partial_\varepsilon \mathcal{J}(x+\varepsilon y)z\big|_{\varepsilon=0}
 =\textstyle\frac12\{\mathcal{A}(z,x)y+\mathcal{A}(z,y)x\}\,.$$
Setting $x=y$ yields $\mathcal{J}(x,x)=\mathcal{J}(x)$. Furthermore
$$
\mathcal{J}(x,y)x=\textstyle\frac12(\mathcal{A}(x,x)y+\mathcal{A}(x,y)x)=-\frac12\mathcal{J}(y)x\,.
$$

Let $A$ be a Jacobi--Tsankov algebraic curvature tensor.
Polarizing the identity
${\mathcal{J}}(x){\mathcal{J}}(y)={\mathcal{J}}(y){\mathcal{J}}(x)$
yields:
$$
\mathcal{J}(x_1,x_2)\mathcal{J}(y_1,y_2)=\mathcal{J}(y_1,y_2)\mathcal{J}(x_1,x_2)\,.
$$
We have $\mathcal{J}(x)x=\mathcal{A}(x,x)x=0$. We prove Assertion (1) by computing:
\begin{eqnarray*}
0=\mathcal{J}(x,y)\mathcal{J}(x,x)x=\mathcal{J}(x,x)\mathcal{J}(x,y)x
 =-\textstyle\frac12\mathcal{J}(x)\mathcal{J}(x)y\,.
\end{eqnarray*}

Since the Jacobi operator is nilpotent, $\{0\}$ is the only eigenvalue of $\mathcal{J}$. This shows that $A$ is
Osserman.

If $p=0$, then ${\mathcal{J}}(x)$ is diagonalizable. Consequently,
${\mathcal{J}}(x)^2=0$ implies ${\mathcal{J}}(x)=0$ for all $x$.
It now follows $A=0$.
 If $p=1$, then $A$ is Osserman implies $A$ has constant sectional
curvature \cite{BBG97,GKV97}. Since ${\mathcal{J}}(x)^2=0$, this
again implies $A=0$. This completes the proof of Theorem
\ref{thm-1.4}. \hfill\qedbox

\medbreak In fact, it is possible to work in a slightly more general
setting. Following Bokan \cite{B90}, one says that
$\mathcal{C}$ is a {\it generalized curvature operator} if it has
the symmetries of the curvature operator defined by a torsion free
connection, i.e. if
\begin{eqnarray*}
&&\mathcal{C}(x,y)z=-\mathcal{C}(y,x)z,\\
&&\mathcal{C}(x,y)z+\mathcal{C}(y,z)x+\mathcal{C}(z,x)y=0.
\end{eqnarray*}
The proof given
above then generalizes immediately to yield:

\begin{corollary}\label{cor-2.1}
If $\mathcal{C}$ is a generalized curvature operator on $V$ which
is Jacobi--Tsankov, then $\mathcal{J}_C$ is Osserman and
$\mathcal{J}_C(x)^2=0$ for all $x\in V$.\end{corollary}

{ Let $\phi$ be a skew-symmetric endomorphism of $V$. Define
$$
A_\phi(x,y,z,w):=\langle\phi y,z\rangle\langle\phi
x,w\rangle-\langle\phi x,z\rangle\langle\phi y,w\rangle
-2\langle\phi x,y\rangle\phi z,w\rangle\,.
$$
The associated Jacobi operator is then given by
$$\mathcal{J}_\phi(x)y=-3\langle
y,\phi x\rangle\phi x\,.$$ In the following example, we exhibit an
algebraic curvature tensor so that $\mathcal{J}(x)^2=0$
for all $x\in V$, but which is not Jacobi--Tsankov. Let
$\mathbb{R}^{(p,q)}$ denote Euclidean space with a metric of
signature $(p,q)$.}

\begin{lemma}\label{lem-2.2}
\ \begin{enumerate} \item There exist skew-symmetric
endomorphisms $\{\phi_1,\phi_2\}$ of $\mathbb{R}^{(4,4)}$ so that
$$\phi_1^2=\phi_2^2=\phi_1\phi_2+\phi_2\phi_2=0,\quad\text{and}\quad\phi_1\phi_2\ne0\,.$$
\item Set $A=-\frac13\{A_{\phi_1}+A_{\phi_2}\}$. Then
$\mathcal{J}_A(x)^2=0$ for all $x$. Furthermore, $A$ is not
Jacobi--Tsankov.
\end{enumerate}
\end{lemma}

\begin{proof}

We apply Lemma 1.4.5 of \cite{G02} to find a collection
$\{e_1,e_2,e_3,e_4\}$ of  skew-symmetric endomorphisms of
$\mathbb{R}^{(4,4)}$ so that:
$$e_1^2=e_2^2=\operatorname{id},\quad
  e_3^2=e_4^2=-\operatorname{id},\quad
  e_ie_j+e_je_i=0\text{ for }i\ne j\,.$$
Set $\phi_1=e_1+e_3$, $\phi_2=e_2+e_4$. These are skew-symmetric endomorphisms with
$$\phi_1^2=\phi_2^2=0,\quad\phi_1\phi_2+\phi_2\phi_1=0\,.$$
Suppose that
$$\alpha:=\phi_1\phi_2=(e_1+e_3)(e_2+e_4)=0\,.$$
We argue for a contradiction. Conjugating by $e_1$ yields
$$e_1\alpha e_1=(-e_1+e_3)(e_2+e_4)=0\,.$$
Adding this equation  to the previous one implies
$e_3(e_2+e_4)=0$. Multiplying by $e_3$ implies $e_2+e_4=0$.
Conjugating this identity by $e_2$ yields $e_2-e_4=0$ and
thus $e_2=0$. This is not possible. Assertion (1) now follows.

To prove Assertion (2), we compute:
\begin{eqnarray*}
&&\mathcal{J}_A(x)y=\langle y,\phi_1x\rangle\phi_1x+\langle y,\phi_2x\rangle\phi_2x,\\
&&\mathcal{J}_A(x_1)\mathcal{J}_A(x_2)y=\langle
y,\phi_1x_2\rangle\langle \phi_1x_2,\phi_1x_1\rangle\phi_1x_1
   +\langle y,\phi_1x_2\rangle\langle \phi_1x_2,\phi_2x_1\rangle\phi_2x_1\\
&&\qquad\qquad+\langle y,\phi_2x_2\rangle\langle
\phi_2x_2,\phi_1x_1\rangle\phi_1x_1
+\langle y,\phi_2x_2\rangle\langle \phi_2x_2,\phi_2x_1\rangle\phi_2x_1\\
&&\qquad\qquad=\langle y,\phi_1x_2\rangle\langle
\phi_1x_2,\phi_2x_1\rangle\phi_2x_1 +\langle
y,\phi_2x_2\rangle\langle \phi_2x_2,\phi_1x_1\rangle\phi_1x_1\,.
\end{eqnarray*}
Since
$$\langle\phi_1x,\phi_2x\rangle=-\langle\phi_2\phi_1x,x\rangle=\langle\phi_1\phi_2x,x\rangle=-\langle\phi_2x,\phi_1x\rangle,$$
we have ${\mathcal{J}}(x){\mathcal{J}}(x)=0$ as desired.

Choose $x_1$ so $\phi_2\phi_1x_1\ne0$. Set $y=\phi_1x_1$. We then
have:
\begin{eqnarray*}
\mathcal{J}_A(x_1)\mathcal{J}_A(x_2)y&=&\langle \phi_1x_1,\phi_1x_2\rangle\langle \phi_1x_2,\phi_2x_1\rangle\phi_2x_1\\
&+&\langle \phi_1x_1,\phi_2x_2\rangle\langle \phi_2x_2,\phi_1x_1\rangle\phi_1x_1\\
&=&\langle \phi_1x_1,\phi_2x_2\rangle^2\phi_1x_1,\\
\mathcal{J}_A(x_2)\mathcal{J}_A(x_1)y&=&\langle \phi_1x_1,\phi_1x_1\rangle\langle \phi_1x_1,\phi_2x_2\rangle\phi_2x_2\\
&+&\langle \phi_1x_1,\phi_2x_1\rangle\langle \phi_2x_1,\phi_1x_2\rangle\phi_1x_2\\
&=&0\,.
\end{eqnarray*}
Choose $x_2$ so $\langle\phi_1x_1,\phi_2x_2\rangle\ne0$. Then
$\mathcal{J}_A(x_1)\mathcal{J}_A(x_2)y\ne0=\mathcal{J}_A(x_2)\mathcal{J}_A(x_1)y$.
\end{proof}

\section{$2$-step Jacobi nilpotent algebraic curvature tensors}\label{sect-3}

Theorem \ref{thm-1.6} will follow from the following
result:

{ \begin{lemma}\label{lem-3.1} Let
$\mathfrak{M}:=(V,\langle\cdot,\cdot\rangle,A)$ be
Jacobi--Tsankov. Suppose that there exist $x,y\in V$ so that
$\mathcal{J}(x)\mathcal{J}(y)\ne0$.
\begin{enumerate}\item There exists $w\in V$ so that
$$
\langle\mathcal{J}(x)\mathcal{J}(y)w,w\rangle=
\langle\mathcal{J}(y)\mathcal{J}(w)x,x\rangle=
\langle\mathcal{J}(w)\mathcal{J}(x)y,y\rangle\ne0\,.$$
\item Let $\mathcal{J}_x:=\mathcal{J}(x)$,
$\mathcal{J}_y:=\mathcal{J}(y)$ and
$\mathcal{J}_{xy}:=\mathcal{J}(x,y)$. Set
$$
\begin{array}{llll}
e_2=\mathcal{J}_x\mathcal{J}_yw,&e_3=\mathcal{J}_xw,&e_4=\mathcal{J}_yw,&e_5=\mathcal{J}_{xy}w\\
f_2=\mathcal{J}_y\mathcal{J}_wx,&f_3=\mathcal{J}_yx,&f_4=\mathcal{J}_wx,&f_5=\mathcal{J}_{yw}x\\
g_2=\mathcal{J}_w\mathcal{J}_xy,&g_3=\mathcal{J}_wy,&g_4=\mathcal{J}_xy,&g_5=\mathcal{J}_{wx}y.
\end{array}
$$
The set $S:=\{w,x,y,e_2,...,e_5,f_2,...,f_5,g_2,...,g_4\}$ is linearly
independent. \item $e_5+f_5+g_5=0$.
 \item $\dim(V)\ge14$.
\end{enumerate}\end{lemma}}

\begin{proof}{ Choose $w$ so that $e_2:=\mathcal{J}(x)\mathcal{J}(y)w\ne0$.
Choose $f$ so $\langle e_2,f\rangle\ne0$.
Set $w(\varepsilon):=w+\varepsilon f$ and
$e_2(\varepsilon):=\mathcal{J}(x)\mathcal{J}(y)w(\varepsilon)$}. Then
$$
p(\varepsilon):=\langle w(\varepsilon),e_2(\varepsilon)\rangle
=\langle w,e_2\rangle+2\varepsilon\langle e_2,f\rangle
+\varepsilon^2\langle \mathcal{J}(x)\mathcal{J}(y)f,f\rangle\,.
$$
As $\langle e_2,f\rangle\ne0$, $p(\varepsilon)$ is a non-trivial
polynomial in $\varepsilon$. Thus it is non-zero for a suitable
choice of $\varepsilon$. Thus { we may choose $w$ so that}
$\langle w,\mathcal{J}(x)\mathcal{J}(y) w\rangle\ne0$. Now,
\begin{eqnarray*}
\langle\mathcal{J}(y)\mathcal{J}(w)x,x\rangle&=&\textstyle
-2\langle\mathcal{J}(y)\mathcal{J}(w,x)w{,x}\rangle=
-2\langle\mathcal{J}(y)w,\mathcal{J}(w,x)x\rangle\\
&=&\langle \mathcal{J}(y)w,\mathcal{J}(x)w\rangle =\langle
\mathcal{J}(x)\mathcal{J}(y)w,w\rangle\,.
\end{eqnarray*}
Similarly, $\langle \mathcal{J}(w)\mathcal{J}(x)y,y\rangle=\langle
\mathcal{J}(x)\mathcal{J}(y)w,w\rangle$ and Assertion (1) follows.

{ Because $\mathcal{J}(x+\varepsilon y)\mathcal{J}(x+\varepsilon
y)=0$ for every $\varepsilon \in \mathbb{R}$ and because
$\mathfrak{M}$ is Jacobi--Tsankov, we have the following
relations:}
$$\begin{array}{lll}
\mathcal{J}_x^2=0,&\mathcal{J}_y^2=0,&
\mathcal{J}_x\mathcal{J}_y=\mathcal{J}_y\mathcal{J}_x,\\
\mathcal{J}_x\mathcal{J}_{xy}=\mathcal{J}_{xy}\mathcal{J}_x=0,&
\mathcal{J}_y\mathcal{J}_{xy}=\mathcal{J}_{xy}\mathcal{J}_y=0,&\mathcal{J}_{xy}^2=-\textstyle\frac12\mathcal{J}_x\mathcal{J}_y\,.
\end{array}$$
We have $\mathcal{J}_w\mathcal{J}_yx\ne0$ and
$\mathcal{J}_w\mathcal{J}_xy\ne0$  by Assertion (1). { To prove
Assertion (2), suppose there is a non-trivial dependence relation
among the elements of $S$:
\begin{eqnarray}
0&=&a_1w+a_2e_2+a_3e_3+a_4e_4+a_5e_5\nonumber\\
&+&b_1x+b_2f_2+b_3f_3+b_4f_4+b_5f_5\nonumber\\
&+&c_1y+c_2g_2+c_3g_3+c_4g_4+c_5g_5\nonumber\\
&=&a_1w+a_2\mathcal{J}_x\mathcal{J}_yw+a_3\mathcal{J}_xw+a_4\mathcal{J}_yw+a_5\mathcal{J}_{xy}w\nonumber\\
&+&b_1x+b_2\mathcal{J}_y\mathcal{J}_wx+b_3\mathcal{J}_yx+b_4\mathcal{J}_wx+b_5\mathcal{J}_{yw}x\label{eqn-3.a}\\
&+&c_1y+c_2\mathcal{J}_w\mathcal{J}_xy+c_3\mathcal{J}_wy+c_4\mathcal{J}_xy+c_5\mathcal{J}_{wx}y\,.\nonumber
\end{eqnarray}
Since we are {\bf not} taking $g_5$, we must set
\begin{equation}\label{eqn-3.b}
c_5=0\,.
\end{equation}}
We can apply $J_xJ_y$ to Equation (\ref{eqn-3.a}) to see
$a_1e_5=0$. Since $e_5\ne0$, $a_1=0$. Similarly $b_1=c_1=0$. If we
now apply $J_x$ to Equation (\ref{eqn-3.a}), we see
\begin{eqnarray*}
&&a_4\mathcal{J}_x\mathcal{J}_yw+c_3\mathcal{J}_x\mathcal{J}_wy=0\quad\text{so}\\
&&0=\langle
a_4\mathcal{J}_x\mathcal{J}_yw+c_3\mathcal{J}_x\mathcal{J}_wy,w\rangle=a_4\langle\mathcal{J}_x\mathcal{J}_yw,w\rangle\,.
\end{eqnarray*}
Since $\langle\mathcal{J}_x\mathcal{J}_yw,w\rangle\ne0$, $a_4=0$.
Similarly, we get $a_3=b_3=b_4=c_3=c_4=0$. Thus Equation
(\ref{eqn-3.a}) simplifies to become
\begin{eqnarray*}
0&=&a_2\mathcal{J}_x\mathcal{J}_yw+a_5\mathcal{J}_{xy}w+b_2\mathcal{J}_y\mathcal{J}_wx+b_5\mathcal{J}_{yw}x
   +c_2\mathcal{J}_w\mathcal{J}_xy+c_5\mathcal{J}_{wx}y\,.
\end{eqnarray*}
Applying $\mathcal{J}_{xy}$ then yields
\begin{eqnarray*}
0&=&a_5\mathcal{J}_{xy}^2w+b_5\mathcal{J}_{xy}\mathcal{J}_{yw}x+c_5\mathcal{J}_{xy}\mathcal{J}_{wx}y\\
&=&(a_5\mathcal{J}_{xy}^2+\textstyle\frac14(b_5+c_5)\mathcal{J}_{x}\mathcal{J}_y)w\\
&=&(a_5-\textstyle\frac12(b_5+c_5))\mathcal{J}_{xy}^2w\,.
\end{eqnarray*}
This shows $a_5=\textstyle\frac12(b_5+c_5)$ or $a_5=b_5=c_5$. { By
Equation (\ref{eqn-3.b}),} we have $a_5=b_5=0$. Taking the inner
product with $x$, $y$, and $w$ then yields, respectively $b_2=0$,
$c_2=0$, and $a_2=0$, which completes the proof of Assertion (2).

{ To prove Assertion (3), we compute:
\begin{eqnarray*}
&&e_5+f_5+g_5=\mathcal{J}_{xy}w+\mathcal{J}_{yw}x+\mathcal{J}_{wx}y\\
&=&\textstyle\frac12\{\mathcal{R}(w,x)y+\mathcal{R}(w,y)x+\mathcal{R}(x,y)w+\mathcal{R}(x,w)y
+\mathcal{R}(y,w)x+\mathcal{R}(y,x)w\}\\
&=&0\,.
\end{eqnarray*}

Assertion (4) is immediate from Assertion (2).}
\end{proof}

The following example in signature $(8,6)$ was motivated by the
proof of Lemma \ref{lem-3.1}. It shows the inequality $\dim(V)\leq
13$ in Theorem \ref{thm-1.6} is sharp. The proof is a computer
assisted { calculation} which we omit in the interest of brevity.
Details are available upon request from the first author.
\begin{lemma}\label{lem-3.2} Let
$\{e_1,\dots,e_4,\bar{e}_1,\dots,\bar{e}_4,\tilde{e}_1,\dots,\tilde{e}_4,f_1,f_2\}$
be a basis for a $14$ dimensional vector space $V$. Relative to this basis, define an inner product
$\langle\cdot,\cdot\rangle$ and an algebraic curvature tensor
$A$ on $V$ whose non-zero components are given up to the usual
$\mathbb{Z}_2$ symmetries by:
\begin{eqnarray*}
&&\textstyle
\langle e_1,e_2\rangle=\langle e_3,e_4\rangle=\langle\bar e_1,\bar e_2\rangle=\langle\bar e_3,\bar e_4\rangle
=\langle\tilde e_1,\tilde e_2\rangle=\langle\tilde e_3,\tilde e_4\rangle=1,\\ 
&&\langle f_1,f_1\rangle=\langle f_2,f_2\rangle=-\textstyle\frac12,\qquad
  \langle f_1,f_2\rangle=\frac14,\\
&&A(e_1,\tilde{e}_1,\tilde{e}_1,e_3)=A(e_1,\bar{e}_1,\bar{e}_1,e_4)=1,\ 
A(\bar{e}_1,e_1,e_1,\bar{e}_3)=A(\bar{e}_1,\tilde{e}_1,\tilde{e}_1,\bar{e}_4)=1,\\
&&A(\tilde{e}_1,e_1,e_1,\tilde{e}_3)=A(\tilde{e}_1,\bar{e}_1,\bar{e}_1,\tilde{e}_4)=1,\\
&&A(e_1,\bar{e}_1,\tilde{e}_1,f_1)=A(e_1,\tilde{e}_1,\bar{e}_1,f_1)=
A(\bar{e}_1,\tilde{e}_1,e_1,f_2)=A(\bar{e}_1,e_1,\tilde{e}_1,f_2)=-\textstyle\frac{1}{2}.
\end{eqnarray*}
Then ${\mathfrak{M}}:=(V,\langle \cdot,\cdot\rangle,A)$ is
Jacobi--Tsankov, { $\mathfrak{M}$ has signature $(8,6)$, and
$\mathfrak{M}$} is not $2$-step Jacobi nilpotent.
\end{lemma}

\section{The classification of indecomposable $2$-step Jacobi nilpotent algebraic curvature tensors}\label{sect-4}

{ In this section, we prove Theorem \ref{thm-1.9}. The following
Lemma shows that Assertion (3) implies Assertion (2) in Theorem
\ref{thm-1.9}.

\begin{lemma}\label{lem-4.1}Let $\mathfrak{M}$ be as in Definition \ref{defn-1.8}. Then $\mathfrak{M}$ is
indecomposable and $2$-step skew-curvature nilpotent.
\end{lemma}

\begin{proof} Suppose there is a non-trivial decomposition $\mathfrak{M}=\mathfrak{M}_1\oplus\mathfrak{M}_2$.
This would then induce a non-trivial decomposition of $(W,A_W)$.
Since $(W,A_W)$ is assumed indecomposable, either $W\subset V_1$
or $W\subset V_2$; we suppose without loss of generality that
$W\subset V_1$. Since $V_2\perp V_1$ and $W\subset V_1$, $V_2\perp W$ so $V_2\subset W^\perp=W$. Thus
$V_2$ is totally isotropic which is false. This shows
$\mathfrak{M}$ is indecomposible. The following argument shows that $\mathfrak{M}$ is $2$-step
curvature nilpotent. Choose a basis
$\{e_i\}$ for
$W$ and choose a basis $\{\bar e_i\}$ for $\bar W$ so the only non-zero components of the inner product
are
$\langle e_i,\bar e_j\rangle=\delta_{ij}$. We have
$$
{\mathcal{A}}(e_i,e_j)e_k=\textstyle\sum_l{
A_W}(e_i,e_j,e_k,e_l)\bar e_l,
$$
while ${\mathcal{A}}(e_i,e_j)e_k=0$ if any entry belongs to $\bar
W$.\end{proof}

{ We now show  Assertion (2) implies Assertion (1) in Theorem
\ref{thm-1.9}.}

\begin{lemma}\label{lem-4.2}
Let $\mathfrak{M}:=(V,\langle\cdot,\cdot\rangle,A)$. If
$\mathfrak{M}$ is $2$-step skew-curvature nilpotent{, }then
$\mathfrak{M}$ is $2$-step Jacobi nilpotent.
\end{lemma}

\medbreak\noindent{\it Proof.}
Suppose $A$ is a $2$-step skew-curvature nilpotent algebraic
curvature tensor. Then $
\mathcal{A}(x_1,x_2)\mathcal{A}(x_3,x_4)=0$ for all
$x_1,x_2,x_3,x_4 \in V$. Hence
\medbreak\qquad $0=-\langle
\mathcal{A}(x_1,x_2)\mathcal{A}(x_3,x_4)x_4,x_2\rangle=-\langle
\mathcal{A}(x_1,x_2)\mathcal{J}(x_4)x_3,x_2\rangle$
\smallbreak\qquad\qquad\qquad$=
    -\langle\mathcal{J}(x_4)x_3,\mathcal{A}(x_1,x_2)x_2\rangle=\langle
\mathcal{J}(x_4)x_3,\mathcal{J}(x_2)x_1\rangle$
\smallbreak\qquad\qquad\qquad$=
\langle\mathcal{J}(x_2)\mathcal{J}(x_4)x_3,x_1\rangle\,.$\hfill\qedbox

\medbreak{ Before completing the proof of Theorem \ref{thm-1.9},
we must establish a technical result.}

\begin{lemma}\label{lem-4.3}
Let $\mathfrak{M}:=(V,\langle\cdot,\cdot\rangle,A)$. { Suppose
that $\mathcal{J}(x)y=0$ for all $x\in V$. Then
$A(x_1,x_2,x_3,y)=0$ for all $x_i\in V$.}\end{lemma}

{ \medbreak\noindent{\it Proof.} We compute:
\begin{eqnarray*}
&&{ A}(x_1,x_2,x_3,y)+{ A}(x_1,x_3,x_2,y)
=2\langle\mathcal{J}(x_2,x_3)x_1,y\rangle\\
&=&2\langle
x_1,\mathcal{J}(x_2,x_3)y\rangle=0\,.
\end{eqnarray*}
Consequently ${ A}(x_1,x_2,x_3,y)=-{ A}(x_1,x_3,x_2,y)$ for all
$x_i\in V$.} This implies
\begin{eqnarray*}
0&=&{ A}(x_1,x_2,x_3,y)+{ A}(x_2,x_3,x_1,y)+{ A}(x_3,x_1,x_2,y)\\
&=&{ A}(x_1,x_2,x_3,y)-{ A}(x_2,x_1,x_3,y)-{ A}(x_1,x_3,x_2,y)\\
&=&{ A}(x_1,x_2,x_3,y)+{ A}(x_1,x_2,x_3,y)+{ A}(x_1,x_2,x_3,y)\\
&=&3{
A}(x_1,x_2,x_3,y)\,.\qquad\qquad\qquad\qquad\qquad\qquad\qquad\qquad\qquad\qedbox
\end{eqnarray*}

{ We complete our discussion by showing that Assertion (1) implies
Assertion (3) in Theorem \ref{thm-1.9}.} Suppose that
$\mathfrak{M}$ is indecomposable and that $\mathfrak{M}$ is
$2$-step Jacobi nilpotent. Set
$$
\bar W:=\operatorname{Span}_{v_1,v_2\in V}\{{\mathcal{J}}(v_1)v_2\}\quad\text{and}\quad
U:=\{v\in V:\mathcal{J}(v_1)v=0\;\forall v_1\in V\}\,.
$$
Then by assumption, $\bar W\subset U$. Furthermore, by Lemma
\ref{lem-4.3}, ${ A}(v_1,v_2,v_3,v_4)=0$ if any of the $v_i\in U$.
Choose a complementary subspace $W_1$ so that $V=U\oplus W_1$.

If $\bar w\in \bar W$, then $\bar
w=\sum_j{\mathcal{J}}(x_j)y_{j}$. Thus if $u\in U$,
\begin{equation}\label{eqn-4.a}
\langle\bar
w,u\rangle=\textstyle\langle\sum_j{\mathcal{J}}(x_j)y_j,u\rangle=\sum_j\langle
y_j,{\mathcal{J}}(x_j)u\rangle=0\,.
\end{equation}
Since the metric is non-degenerate, there must exist ${ \tilde
w}\in W_1$ so $\langle{\tilde w},\bar w\rangle\ne0$. Thus the
natural map $W_1\rightarrow \bar W^*$ defined by
$\langle\cdot,\cdot\rangle$ is surjective. Let $\{\bar
w_1,...,\bar w_k\}$ be a basis for $\bar W$. Choose elements
$\{\tilde w_1,...,\tilde w_k\}$ in $W_1$ so
$$\langle\tilde w_i,\bar w_j\rangle=\delta_{ij}\,.$$

Suppose that $\{\tilde w_1,...,\tilde w_k\}$ do not span $W_1$. We
may then choose $0\ne\tilde w\in W_1$ so that $\tilde w\perp \bar
W$. Since $\tilde w\notin U$, there exists $y$ so that
${\mathcal{J}}(y)\tilde w\ne0$. Choose $z\in V$ so
$$0\ne\langle {\mathcal{J}}(y)\tilde w,z\rangle=\langle\tilde w,{\mathcal{J}}(y)z\rangle\,.$$
This contradicts the fact that $\tilde w\perp \bar W$.
Thus $\{\tilde w_1,...,\tilde w_k\}$ is a basis for $W_1$. We set
$$
w_i:=\tilde w_i-\textstyle\frac12\sum_j\langle\tilde
w_i,\tilde w_j\rangle\bar w_j\quad\text{and}\quad
W:=\operatorname{Span}\{w_i\}\,.
$$
Then the relations of Equation (\ref{eqn-1.a}) are satisfied. Furthermore,
$$V=U\oplus W\,.$$

Let $\{\bar w_1,...,\bar w_k,\tilde u_1,...,\tilde u_l\}$ be a
basis for $U$. By Equation (\ref{eqn-4.a}), $\langle\bar
w_i,{\tilde u_j}\rangle=0$. Set
$$u_i:=\tilde u_i-\textstyle\sum_j\langle w_j{,}\tilde u_i\rangle\bar w_j\,.$$
We then have $\langle u_i,w_i\rangle=\langle u_i,\bar
w_i\rangle=0$. Let $T:=\operatorname{Span}\{u_i\}$. { Then:}
$$ (V,\langle\cdot,\cdot\rangle,{ A}) =(W\oplus\bar
W,\langle\cdot,\cdot\rangle|_{W\oplus\bar W}, A|_W\oplus 0)\
\oplus\ (T,\langle\cdot,\cdot\rangle|_T,0)\,.
$$
Since $(V,\langle\cdot,\cdot\rangle,{ A})$ is indecomposable,
$T=\{0\}$.\hfill\qedbox

\section*{Acknowledgments}
The research of M. Brozos-V\'azquez was partially supported by project BFM 2003-02949 (Spain).  The research of P. Gilkey was
partially supported by the Max Planck Institute for the Mathematical Sciences (Leipzig, Germany).  It is a
pleasure for both authors to acknowledge helpful conversations with C. Dunn and E. Puffini.

\end{document}